\documentclass[12pt]{article}
\usepackage{fancyhdr, array,calc,graphicx,url,tabularx}
\usepackage{geometry}
\usepackage{latexsym}
\usepackage{amssymb}
\usepackage{amsmath}
\usepackage{makeidx}
\usepackage{undertilde}

\makeindex

{
   \end{minipage}
   \vspace*{\stretch{3}}
   \clearpage
}

\renewcommand{\gg}{\gamma}

\newcommand{\gG}{\Gamma}

\newcommand{\bR}{{\mathbb{R}}}

\newcommand{\rest}{\restriction}
\newcommand{\la}{\langle}
\newcommand{\ra}{\rangle}

\newcommand{\forces}{\Vdash}

\renewcommand{\models}{\vDash}
\newcommand{\powerset}{{\wp}}
\newcommand{\cp}{{\rm cp }}

\newtheorem{theorem}{Theorem}[section]
\newtheorem{proposition}[theorem]{Proposition}

\newtheorem{lemma}[theorem]{Lemma}
\newtheorem{corollary}[theorem]{Corollary}

\numberwithin{figure}{section}

\newenvironment{proof}{{\it{
Proof.}}}{\nopagebreak\mbox{}{\hfill$\square$}
\par\bigskip}
\newenvironment{Claim}{\noindent{\it
Claim} {}}{}
\newenvironment{Case}{\noindent{\it
Case}}{}

\newcommand{\rprop}[1]{Proposition~\ref{#1}}
\newcommand{\rthm}[1]{Theorem~\ref{#1}}
\newcommand{\rlem}[1]{Lemma~\ref{#1}}

\newcommand{\rcor}[1]{Corollary~\ref{#1}}

\def\inseg{\trianglelefteq}

\def\k{\kappa}
\def\a{\alpha}

\def\d{\delta}

\def\l{\lambda}

\def\P{{\mathcal{P} }}
\def\W{{\mathcal{W} }}
\def\Q{{\mathcal{ Q}}}

\def\R{{\mathcal R}}

\def\M{{\mathcal{M}}}
\def\N{{\mathcal{N}}}
\def\F{{\mathcal{F}}}
\def\T {{\mathcal{T}}}
\def\U{{\mathcal{U}}}
\def\S{{\mathcal{S}}}
\def\F{{\mathcal{F}}}

\def\VT{{\vec{\mathcal{T}}}}

\def\cp #1{{ crit  #1 }}

\def\and{\mathrel{\kern1pt\&\kern1pt}}

\def\inseg{\triangleleft}
\def\insegeq{\trianglelefteq}

\def\<#1>{\langle\,#1\,\rangle}

 \input xy
 \xyoption{all}

\begin{document}
\title{The mouse set conjecture for sets of reals
 }
\date{}
\author{Grigor Sargsyan\thanks{First author's work is partially based upon work supported by the National Science Foundation under Grant No DMS-0902628.}\\
        Department of Mathematics\\
        Rutgers University\\
        110 Frelinghuysen Rd.\\
        Piscataway, NJ 08854\\
        http://math.rutgers.edu/$\sim$gs481\\
        grigor@math.rutgers.edu,\\\\
John Steel\\
Department of Mathematics\\
University of California\\
Berkeley, California, 94720 USA\\
http://math.berkeley.edu/$\sim$steel\\
steel@math.berkeley.edu}
\maketitle
\begin{abstract}
We show that the \textit{Mouse Set Conjecture} for sets of reals is true in the minimal model of $AD_{\mathbb{R}}+``\Theta$ is regular". As a consequence, we get that below $AD_{\mathbb{R}}+``\Theta$ is regular", models of $AD^++\neg AD_{\mathbb{R}}$ are hybrid mice over $\mathbb{R}$. Such a representation of models of $AD^+$ is important in core model induction applications.
\end{abstract}

\thispagestyle{empty}

One of the central open problems in descriptive inner model theory is the conjecture known as the \textit{Mouse Set Conjecture} ($MSC$). It conjectures that under $AD^+$ ordinal definable reals are exactly those that appear in $\omega_1$-iterable mice. The counterpart of this conjecture for sets of reals conjectures that under $AD^+$, the sets of reals which are ordinal definable from a real are exactly those that appear in countably iterable mice over $\mathbb{R}$. In \cite{ATHM}, the first author proved that $MSC$ holds in the minimal model of $AD_{\mathbb{R}}+``\Theta$ is regular", but $MSC$ for sets of reals was left open. The goal of this paper is to establish that $MSC$ for sets of reals holds in the minimal model of $AD_{\mathbb{R}}+``\Theta$ is regular". 

We will establish a stronger form of $MSC$ known as the \textit{Strong Mouse Set Conjecture} ($SMSC$). We say $\M$ is countably $\k$-iterable if all of its sufficiently elementary countable substructures are $\k$-iterable. We say $\M$ is countably iterable if $\M$ is countably $\omega_1$-iterable. Thus, under $AD$, if $\M$ is countably iterable then $\M$ is countably $\omega_1+1$-iterable. 

In what follows, we will let ``hod pair" stand for a hod pair below $AD_{\mathbb{R}}+``\Theta$ is regular", i.e., the corresponding hod mouse cannot have inaccessible limit of Woodin cardinals (see Definition 1.34 of \cite{ATHM}). Given an iteration strategy $\Sigma$ for a countable structure, we let $Code(\Sigma)$ be the set of reals coding $\Sigma$ for trees of length $\omega_1$. Given a hod pair $(\P, \Sigma)$ we let 
\begin{center}
$Lp^{\Sigma}(\mathbb{R})=\cup\{ \M: \M$ is a sound countably iterable $\Sigma$-mouse over $\mathbb{R}$ projecting to $\mathbb{R}\}$.
\end{center}
The following is the statement of $SMSC$ for sets of reals. Recall the notions of branch condensation and fullness preservation from \cite{ATHM} (see Definition 2.14 and Definition 2.27 of \cite{ATHM}). Recall that $OD_X$ stands for the class of sets ordinal definable from a finite sequence consisting of members of $X$. \\  

\textbf{The Strong Mouse Set Conjecture for sets of reals, $SMSC(\mathbb{R})$:} Assume $AD^+$. Suppose $(\P, \Sigma)$ is a hod pair such that $\Sigma$ has branch condensation and is fullness preserving. Then
\begin{center}
$\{ A\subseteq \mathbb{R}: \exists x\in \bR (A$ is $OD_{\{\Sigma, x\}})\}=Lp^\Sigma(\bR)$. 
\end{center}

To following is the main theorem of this paper.

\begin{theorem}\label{main theorem}Assume $AD^++V=L(\powerset(\mathbb{R}))$.  Suppose $(\P, \Sigma)$ is a hod pair such that the following holds.
\begin{enumerate}
\item $\P$ does not have inaccessible limit of Woodin cardinals. 
\item $\Sigma$ has branch condensation and is fullness preserving.
\item $MSC$ for $\Sigma$ holds, i.e., for every $x, y\in \mathbb{R}$, $x\in OD(\Sigma, y)$ iff $x$ is in a $\Sigma$-mouse over $y$.
\item Every set of reals $A$ is $OD(\Sigma, x)$ for some real $x$.
\end{enumerate}
Then 
\begin{center}
$\powerset(\mathbb{R})=\powerset(\bR)\cap Lp^\Sigma(\mathbb{R})$.
\end{center}
In particular, $V=L(Lp^\Sigma(\mathbb{R}))$.
\end{theorem}

\begin{corollary}\label{main corollary} Suppose $V=L(\powerset(\bR))$ and $AD^+$ holds. Suppose further that for any $\a$ such that $\theta_\a<\Theta$, letting $\Gamma=\{ A\subseteq \bR: w(A)<\theta_\a\}$, $L(\Gamma, \bR)\models \neg AD_\bR$. Then $SMSC(\bR)$ holds.  In particular, $SMSC(\mathbb{R})$ is true in the minimal model of $AD_{\mathbb{R}}+``\Theta$ is regular". 
\end{corollary}
\begin{proof}
It is shown in \cite{ATHM} that if $(\P, \Sigma)$ is as in the hypothesis of \rthm{main theorem} then clause 3 holds in $L(\Gamma_{\a+1})$ where $\a$ is such that $\theta_\a=w(Code(\Sigma))$ and $\Gamma_{\a+1}=\{ A\subseteq \bR: w(A)<\theta_{\a+1}\}$. It then follows from \rthm{main theorem} that  $\Gamma_\a=\powerset(\bR)\cap Lp^\Sigma(\mathbb{R})$ implying that $SMSC(\bR)$ holds.
\end{proof}

All the background material that we will need in this paper is spelled out in \cite{ATHM}. We assume that our reader is familiar with some aspects of  it. One important comment is that in general hybrid mice over $\mathbb{R}$ or any non-self-wellordered set are not defined (recall that a set $X$ is self-wellordered if there is a wellordering of it in $\mathcal{J}_\omega(X)$). Given an iteration strategy $\Sigma$ with hull condensation, the $\Sigma$-mice over self-wellordered sets are defined according to the following principle. At a typical stage where we would like to add more of $\Sigma$ to the model, we choose the least tree $\T$ for which $\Sigma(\T)$ hasn't been defined. However, $\mathbb{R}$ isn't self-wellordered and hence, we cannot choose the least such $\T$.

 In \cite{ATHM}, the first author gave a definition of premice over any non-self-wellordered sets under the hypothesis that $\M_1^{\#, \Sigma}$ exists, i.e., there is a minimal active $\Sigma$-mouse with one Woodin cardinal (see Definition 3.37 of \cite{ATHM}). This extra assumption is benign as under $AD^+$ whenever $(\P, \Sigma)$ is a hod pair such that $\Sigma$ has branch condensation and is fullness preserving, $\M_1^{\#, \Sigma}$ exists and is $\Theta$-iterable. The proof is the same as the proof that shows that $AD^{L(\bR)}$ implies that $\M_1^{\#}$ exists and is $\Theta$-iterable in $L(\bR)$ (see \cite{OIMT}). One consequence of the indexing of the strategy introduced in Definition 3.37 of \cite{ATHM} is that it allows us to perform $S$-constructions, which we will use in this paper (see Chapter 3 of \cite{ATHM}).

\rcor{main corollary} has been used in core model induction applications. See, for instance, \cite{CuBF}, \cite{thesis}, \cite{UBH} or Chapter 7 of \cite{CMI}.  Before we begin the proof of \rthm{main theorem}, we introduce Prikry tree forcing associated with Martin's measure on degrees. 

\textbf{Acknowledgments:} The author's would like to express their gratitude to the referee for a very helpful list of comments. These comments made the paper significantly better. 

\section{Prikry tree forcing on degrees}

We develop the notion of Prikry forcing that we need in a general context. Assume $ZF-Replacement+AD$. Let $\mathcal{D}$ be the set of Turing degrees. Let $f:\mathcal{D}^{<\omega}\rightarrow HC$ be some function. We would like to define \textit{Prikry tree forcing} on degrees associated to $f$. Let $\mu$ be Martin's measure. We let  $(p,A) \in \mathbb{P}^f$ if 
\begin{enumerate}
\item $p\in \mathcal{D}^{<\omega}$,
\item for any $n<lh(p)$, $(f(p\rest n), p\rest n) \in L[p(n)]$,
\item $A\subseteq \cup_{n<\omega}\mathcal{D}^{<n}$ is a tree with stem $p$ such that for every $q\in A$ (in particular, $p\subseteq q$), 
\begin{center}
$\{ d: q^\frown d\in A\} \in \mu$. 
\end{center}
\end{enumerate}
Given $(p, A), (q, B)\in \mathbb{P}^f$ we let
\begin{center}
$(p,A) \preccurlyeq (q, B)$ iff $p$ end-extends $q$, $A \subseteq B$ and $p\in B$. 
\end{center}
We say $p\in \mathcal{D}^{<\omega}$ is a \textit{precondition} if it satisfies 1 and 2 above. Given a precondition $p$ and $d\in \mathcal{D}$, we say $d$ is \textit{valid} at $p$ if $p^\frown d$ is a precondition.
 
Given a $\mathbb{P}^f$-generic $G$ we let $g=\cup \{p : \exists X(p, X)\in G\}$. We then let
\begin{center}
$G^i=_{def}f(g\rest i+1)$ and $f(G)=_{def}\cup_{i<\omega}G^i$ 
\end{center}

The following is proved by a standard fusion argument.  

\begin{lemma}
\label{Prikry property}
 $\mathbb{P}^f$ has the Prikry property. More precisely,  suppose $Z$ is a countable set of $\mathbb{P}^f$-terms, $\phi$ is a formula, and $(p,A) \in \mathbb{P}^f$. Then there is a condition $(p, W) \in \mathbb{P}^f$ deciding $\phi [\tau]$ for all $\tau \in Z$ such that $W\in OD_{Z, \{f,p, A\}}$.
\end{lemma}
\begin{proof}
 We will show that there is a condition $(p, T_\tau)$ deciding $\phi[\tau]$ such that $\la T_\tau : \tau\in Z\ra\in OD_{Z, \{f,p, A\}}$. It then follows that $(p, \cap_{\tau\in Z}T_\tau)$ is as desired. We say $q$ is \textit{positive} if $(\exists Y) \ ((q,Y) \Vdash \phi[\tau])$, \textit{negative} if $(\exists Y) ((q,Y) \Vdash \neg\phi[\tau])$, and \textit{ambiguous} if it is neither positive nor negative. Notice that $q$ cannot be both positive and negative. Fixing $\tau$, we shrink $A$ to some tree $T$ such that given any $r\in T$ and any one step extensions $q_1, 1_2\in T$ of $r$, both $q_1$ and $q_2$ are simultaneously ambiguous, positive or negative. 
 
  We define a sequence of functions $\la H^i: i<\omega\ra$ such that 
\begin{center}
$dom(H^i)=\{ q: p\trianglelefteq q$ and $q$ is a precondition$\} $ 
\end{center}  
and $rng(H^i)\subseteq \{0,1,2\}$. First define $H$ on $\{ (q, d) : q^\frown d$ is a precondition$\}$ by 
\begin{center}$H(q, d)=\begin{cases}
0 :& q^\frown d\  \text{is positive}\\
1 :& q^\frown d\  \text{is negative}\\
2 :& q^\frown d\  \text{is ambiguous}.
\end{cases}$
\end{center}
Now, let $H^0(q)=i$ if for $\mu$-a.e. $d$ is such that $H(q, d)=i$. Given $\la H^i: i\leq k\ra$ define $H^{k+1}$ by setting
$H^{k+1}(q)=i$ if for $\mu$-a.e. $d$ is such that $H^k(q^\frown d)=i$.

We then define a decreasing sequence of conditions $(p, T^i)$ by induction as follows. We will have that $(p, T^0)\preceq (p, A)$.  We define $T^0$ by induction on the length of conditions. We let $T^0\rest m$ be $T^0$ restricted to sequences of length $m$. Suppose we have defined $T^0\rest m+1$ for $m+1\geq lh(p)$. Given $q\in T^0\rest m$ such that $lh(q)=m$ we let 
\begin{center}
$\{ q^\frown d \in A : H(q, d)=H^0(q)\}$
\end{center}
be the one step extensions of $q$ in $T^0$. This finishes our description of $T^0$. 

Suppose now we have defined $\la (p, T^i) : i\leq k\ra$ and $T^{k+1}\rest m+1$. Given $q\in T^{k+1}\rest m$ such that $lh(q)=m$, we let 
\begin{center}
$\{ q^\frown d \in T^k : H^{k}(q^\frown d)=H^{k+1}(q)\}$
\end{center}
be the one step extensions of $q$ in $T^{k+1}$. This finishes our description of $\la (p, T^i): i\leq \omega\ra$. Let $T_\tau=\cap_{i<\omega} T^i$. 

We claim that $(p, T_\tau)$ decides $\tau$. Suppose not. We then have two conditions $(q, X)$ and $(r, Y)$ such that both are below $(p, T_\tau)$ and
\begin{enumerate}
\item $lh(q)=lh(r)$,
\item $(q, X)\forces \phi[\tau]$,
\item $(r, Y)\forces \neg \phi[\tau]$.
\end{enumerate}
Let now $s$ be the common initial segment of $q$ and $r$. Let $s=(d_i: i\leq m)$, $q=s^\frown (q_i: i<n)$ and $r=s^\frown(r_i:i< n)$. It follows from our construction that 
\begin{center}
$H(s^\frown (q_i: i <n-1), q_{n-1})=H^0(s^\frown (q_i : i\leq n-1))=H^1(s^\frown (q_i: i<n-2))=...=H^{n-1}(s)$\\
$H(s^\frown (r_i: i <n-1), r_{n-1})=H^0(s^\frown (r_i : i\leq n-1))=H^1(s^\frown (r_i: i<n-2))=...=H^{n-1}(s)$.
\end{center}
It then follows that $H(s^\frown (q_i: i < n-1), q_{n-1} )=H(s^\frown (r_i: i < n-1), r_{n-1})$, which is a contradiction.
\end{proof}

We now turn to proving \rthm{main theorem}.

\section{The proof}

We assume $AD^++V=L(\powerset(\mathbb{R}))$ and let $(\P, \Sigma)$ be as in the hypothesis of \rthm{main theorem}. Given a good pointclass $\Gamma$\footnote{i.e., a point class closed under $\exists^\bR$, continuous preimages and images, and having the scale property} and $a\in HC$, we let $Lp^{\Gamma, \Sigma}(a)$ be the union of sound $\Sigma$-mice over $a$ projecting to $a$ whose iteration strategy is coded by a set in $\Gamma$. Our first lemma is an easy lemma. Below, $MC(\Sigma)$ (mouse capturing relative to $\Sigma$) is the statement that for every $x, y\in \bR$, $x\in OD_{\Sigma, y}$ if and only if there is an $\omega_1$-iterable sound $\Sigma$-mouse $\M$ over $y$ such that $x\in \M$.

\begin{lemma}\label{easy lemma}
For any good pointclass $\Gamma\not =\Sigma^2_1(Code(\Sigma))$ there is a good pointclass $\Gamma_1\not =\Sigma^2_1(Code(\Sigma))$ such that $\Gamma \cup\{ Code(\Sigma)\} \subseteq \Gamma_1$ and for any $a\in HC$,
\begin{center}
$C_{\Gamma_1}(a)=Lp^{\Gamma_1, \Sigma}(a)$.
\end{center}
\end{lemma}
\begin{proof} Fix a good pointclass $\Gamma\not =\Sigma^2_1(Code(\Sigma))$. Because $MC(\Sigma)$ holds, using $\Sigma_1(Code(\Sigma))$-reflection, we can find $\Gamma_1$ and $\a$ such that $\Gamma_1\not =\Sigma^2_1(Code(\Sigma))$, $\Gamma \cup\{ Code(\Sigma)\} \subseteq \Delta_{\Gamma_1}$, $\mathcal{J}_\a(\Gamma_1, \bR)\models ZF-Replacement$, $\Gamma_1=(\Sigma^2_1(Code(\Sigma)))^{\mathcal{J}_\a(\Gamma_1, \mathbb{R})}$ and 
\begin{center}
$\mathcal{J}_\a(\Gamma_1, \mathbb{R})\models MC(\Sigma)$.
\end{center}
\end{proof}

%Next we prove that $\utilde{\Delta}^2_1(Code(\Sigma))$ is captured by countably iterable $\Sigma$-mice over reals.
%
%\begin{lemma}\label{capturing below delta21} $\utilde{\Delta}^2_1(Code(\Sigma))\subseteq Lp^\Sigma(\mathbb{R})$ 
%\end{lemma}
%\begin{proof} Fix $A\in \utilde{\Delta}^2_1(Code(\Sigma))$. Let $\l=\utilde{\d}^2_1(Code(\Sigma))$. We have that the Solovay measure on $\powerset_{\omega_1}(\mathbb{R})$ restricted to sets coded by a set in $\utilde{\Delta}^2_1(\Sigma)$ is a normal fine ultrafilter.  Let $\Gamma\not =\Sigma^2_1(\Sigma)$ be a good pointclass such that 
%\begin{center}
%$C_{\Gamma}(a)=Lp^{\Gamma, \Sigma}(a)$.
%\end{center}
%for all $a\in HC$ and $A\in \utilde{\Delta}_\Gamma$. We now have that for $\mu$ many $\sigma$, $A\cap \sigma\in Lp^{\Gamma, \Sigma}(\sigma)$. Let $\M(\sigma)=Lp^{\Gamma, \Sigma}(\sigma)$. We have that $\powerset(\sigma)\cap L(\M(\sigma))=\powerset(\sigma)\cap \M(\sigma)$. Let now $\W$ be the ultrapower of $L(\M(\sigma))$ by $\mu$ where we only use functions coded by 
%a set in $\utilde{\Delta}^2_1(\Sigma)$. Then $\W$ is a $\Sigma$-premouse over $\mathbb{R}$. Clearly, $A\in \W$. To see that $\W$ is countably iterable notice that for any $\pi: \Q\rightarrow \W|\eta$ such that $\Q$ is countable, there are $\mu$ many $\sigma\in \powerset_{\omega_1}(\mathbb{R})$ and ordinals $\eta_\sigma$ such that there is $\tau:\Q\rightarrow L(\M(\sigma))|\eta_\sigma$. It then follows that $\Q$ is $(\omega, \omega_1)$-iterable.
%\end{proof}

Suppose now that $\powerset(\mathbb{R})\not =Lp^\Sigma(\mathbb{R})$. Using $\Sigma_1$-reflection we get $\Gamma\subset\utilde{\Delta}^2_1(Code(\Sigma))$ and $\a<\utilde{\delta}^2_1(Code(\Sigma))$ such that
\begin{enumerate}
\item $\Gamma=\powerset(\mathbb{R})\cap \mathcal{J}_\a(\Gamma, \mathbb{R})$ and $\a$ ends a $\Sigma_1$-gap,
\item $\mathcal{J}_\a(\Gamma, \mathbb{R})\models \phi$ where $\phi$ is the conjunction of the following statements:
\begin{enumerate}
\item $ZF-Replacement + DC_\mathbb{R}+MC(\Sigma)$,
%\item $\exists U (U \textrm{ is a } \Sigma^2_1 \textrm{ universal set})$
\item there is an $OD$ set of reals $A$ such that $A\not \in Lp^\Sigma(\mathbb{R})$.
\end{enumerate}
\end{enumerate}

We let $N=\mathcal{J}_\a(\Gamma, \mathbb{R})$. Let $U$ be the set of pairs $(x, y)\in \bR^2$ such that  $y$ codes a sound $\Sigma$-mouse $\M$ over $x$ that projects to $x$ and has an $\omega_1$-iteration strategy in $N$.

 Since $MC(\Sigma)$ holds in $N$, $U$ is a universal $(\Sigma^2_1(Code(\Sigma)))^N$-set. Let $A\in N$ be an $OD$ set of reals witnessing clause (b) of $\phi$. We assume that $A$ has the minimal Wadge rank among the sets witnessing clause b of $\phi$. Using the results of Chapter 3 of \cite{JacksonHST}, we can get $\vec{B}=\la B_i: i<\omega\ra$ which is a semiscale on $U^c$ such that each $B_i\in (OD_{\Sigma})^N$. The following fact is a well known consequence of $MC(\Sigma)$. 
 
 \begin{proposition}\label{mice beyond n} The following statements are true.
 \begin{enumerate}
 \item There is a cone of $x$ such that there is $\M\trianglelefteq Lp^\Sigma(x)$ such that $\rho_\omega(\M)=x$ and $\M$ doesn't have an iteration strategy in $N$.
 \item Let $x$ be a base of the above cone. Then for every $a\in HC$ such that $x\in \mathcal{J}_{\omega}(a)$, there is $\M\trianglelefteq Lp^\Sigma(a)$ such that $\rho(\M)=a$ and $\M$ doesn't have an iteration strategy in $N$.
 \end{enumerate}
 \end{proposition}
 \begin{proof}
 Clause 2 follows from clause 1. To see this, fix a real $x$ such that it is base for the cone of clause 1. Then whenever $a\in HC$ is such that $x\in a$ and $y$ is a real coding $a$ generically over $Lp^\Sigma(a)$, then $Lp^{\Sigma}(a)[y]=Lp^\Sigma(y)$. Indeed, this follows from $S$-constructions (see Section 2.11 of \cite{ATHM}).
 
Clause 1 is an easy consequence of $MC(\Sigma)$. Indeed, suppose clause 1 fails. Then\\\\
(1) there is an $x\in \bR$ such that for all $y\in \bR$ such that $x\leq_T y$, $Lp^{\Sigma}(y)=(Lp^{\Sigma}(y))^N$.\\\\
Because $MC(\Sigma)$ holds, letting $C$ be the set of pairs $(y, z)\in \bR^2$ such that $x\leq_T y$ and $z$ codes an $\omega_1$-iterable sound $\Sigma$-mouse $\M$ over $y$ projecting to $y$, $C$ is a universal $\Sigma^2_1(Code(\Sigma))$ set. Because $\Gamma\subset\utilde{\Delta}^2_1(Code(\Sigma))$, we cannot have that $C$ is the universal $(\Sigma^2_1(Code(\Sigma)))^N$ set. It follows from (1), however, that $C\in N$ and $N\models ``C$ is the universal $\Sigma^2_1(Code(\Sigma))$-set", contradiction.
 \end{proof}

Let now $x$ be a base of the cone from clause 1 of \rprop{mice beyond n}. We say $a$ is \textit{good} if $a\in HC$ and $x\in \mathcal{J}_{\omega}(a)$. For each good $a$ let $\M(a)$ be the least $\Sigma$-mouse with no iteration strategy in $N$. Let $F$ be the set of pairs $(a, \M(a))$. It follows that if $F^*$ is the set of reals coding $F$ then $F^*\in \utilde{\Delta}^2_1(Code(\Sigma))$. Furthermore, there is a set $C\in \utilde{\Delta}^2_1(Code(\Sigma))$ such that for every good $a$, the set of reals coding the unique iteration strategy of $\M(a)$ is Wadge reducible to $C$. Let then 
\begin{center}
$D=\{(y, \sigma)\in \mathbb{R}^2: y$ codes a good $a$ and $\sigma$ codes a continuous function $f$ such that $f^{-1}[C]$ is the iteration strategy of $\M(a)\}$.
\end{center}
 We have that $D\in \utilde{\Delta}^2_1(Code(\Sigma))$.

 Let $\Gamma_1$ be a good pointclass such that $F^*, Code(\Sigma), \vec{B}, U, C, D \in \utilde{\Delta}_{\Gamma_1}$. Moreover, it follows from \rlem{easy lemma} that we can require that for any $a\in HC$
 \begin{center}
 $C_{\Gamma_1}(a)=Lp^{\Gamma_1, \Sigma}(a)$
 \end{center}
 
 Let now $(\N^*_z, \d_z,\Sigma_z)$ be as in Theorem 1.2.9 of \cite{ATHM} with the property that $(\N^*_z, \d_z,\Sigma_z)$ Suslin, co-Suslin captures $Code(\Sigma), \vec{B}, U, C, D$ (where Suslin capturing is defined on page 36 of \cite{ATHM}, also see the next paragraph) \footnote{For convenience, we restate Theorem 1.2.9 of \cite{ATHM}. \begin{theorem}[Woodin, Theorem 10.3 of \cite{DMATM}]\label{n*x} Assume $ AD^+$. Suppose $\Gamma$ is a good pointclasses and there is a good pointclass $\gG^*$ such that $\Gamma\subseteq \Delta_{\gG^*}$. Suppose $(N, \Psi)$ Suslin, co-Suslin capture $\Gamma$. There is then a function $F$ defined on $\mathbb{R}$ such that for a Turing cone of $x$, $F(x)=(\N^*_x, \M_x, \d_x, \Sigma_x)$ such that
\begin{enumerate}
\item $N\in L_1[x]$,
\item $\N^*_x|\d_x=\M_x| \d_x$,
\item $\M_x$ is a $\Psi$-mouse: in fact, $\M_x=\M_1^{\Psi, \#}(x)|\k_x$ where $\k_x$ is the least inaccessible cardinal of $\M_1^{\Psi, \#}$.\footnote{$\M_1^{\Psi, \#}$ is the minimal $\Psi$-mouse having a Woodin cardinal and a last extender.},
\item $\N^*_x\models ``\d_x$ is the only Woodin cardinal",
\item $\Sigma_x$ is the unique iteration strategy of $\M_x$,
\item $\N^*_x= L(\M_x, \Lambda)$ where $\Lambda$ is the restriction of $\Sigma_x$ to stacks $\VT\in \M_x$ that have finite length and are based on $\M_x\rest \d_x$,
\item $(\N^*_x, \Sigma_x)$ Suslin, co-Suslin captures $Code(\Psi)$ and hence, $(\N^*_x, \Sigma_x)$ Suslin, co-Suslin captures $\Gamma$,
\item $(\N^*_x, \d_x, \Sigma_x)$ is a self-capturing background triple. 
\end{enumerate}
\end{theorem}}. We have that for any $\eta<\d_z$, $C_{\Gamma_1}(\N^*_z|\eta)\in \N^*_z$. 
 
Let $\Phi=(\utilde{\Sigma}^2_1)^N$. We have that $\Phi$ is a good pointclass. Because $\vec{B}$ is Suslin captured by $\N^*_z$, we have $(\delta_z^+)^{\N^*_z}$-complementing trees $T, S\in \N^*_z$ which capture $\vec{B}$ in the sense that whenever $i: \N^*_z\rightarrow \N$ is an iteration embedding according to $\Sigma_z$ and $g$ is a generic over $\N$ for a poset of size $\leq i(\d_z)$ then 
\begin{center}
$p[i(T)]\cap \N[g]= B\cap \N[g]$.
\end{center}
Let $\k$ be the least cardinal of $\N^*_z$ which is $<\d_z$-strong in $\N^*_z$. 
 
 Next, we fix a notation. For each $a\in HC$, we let $\W(a)=Lp^{\Gamma, \Sigma}(a)$. Using the results of Section 2.11 of \cite{ATHM}, we have that if $g\subseteq Coll(\omega, a)$ is $\W(a)$-generic then 
\begin{center}
$\W(a)[g]=\W(a, g)=\W(x_g)$\ \ \ \ \ \ \ \ \ \ \ (1).
\end{center}
where $x_g$ is the generic real coding $a$. The following claim is standard.\\

\begin{lemma}\label{many woodin cardinal}
 $\N^*_z\models ``\k$ is a limit of cardinals $\eta$ such that $\eta$ is  a Woodin cardinal in $\W(\N^*_z|\eta)$".
 \end{lemma}
\begin{proof}
Working in $\N^*_z$, let $\lambda=\d_z^{++}$ and let $\pi: M\rightarrow \N^*_z|\l$ be an elementary substructure such that
\begin{enumerate}
\item $T,S\in ran(\pi)$ and 
\item letting $\cp(\pi)=\eta$, $V_\eta^{\N^*_z}\subseteq M$, $\pi(\eta)=\d_z$ and $\eta>\kappa$.
\end{enumerate}
By elementarity, we have that $M\models ``\eta$ is Woodin". Letting $\pi^{-1}(\la T, S\ra)=\la \bar{T}, \bar{S}\ra$, we have that $(\bar{T},\bar{S})$ Suslin captures $\vec{B}$ over $M$ at $(\eta^+)^M$. This implies that whenever $a\in M|(\eta^+)^M$, $\W(a)\in M$. To see this, first note that we have that whenever $g\subseteq Coll(\omega, a)$ is $M$-generic and $x_g$ is the generic real then $\W(x_g)\in M$. But using (1) above, we have that $\W(x_g)=\W(a)[g]$. Therefore, $\W(a)\in M[g]$. Since $g$ was arbitrary, we have that $\W(a)\in M$.

We now have that $\W(\N^*|\eta)\in M$ and since $M\models ``\eta$ is Woodin", we have that $\W(\N^*|\eta)\models ``\eta$ is Woodin".  Because $\k$ is $<\d_z$-strong in $\N^*_z$ and because $a\rightarrow \W(a)$ is definable over $\N^*_z$, we have that for unboundedly many $\nu<\eta$, $\W(\N^*_z|\nu)\models ``\nu$ is Woodin". 
\end{proof}

\subsection{A $\Sigma$-mouse beyond $N$}

In this section, we prove the following important lemma.

\begin{lemma}\label{a mouse beyond n}
There is a $\Sigma$-mouse $\N$ such that there is a sequence $(\gamma_i : i<\omega)$ with the property that
\begin{enumerate}
\item $(\gamma_i : i<\omega)$ is the sequence of Woodin cardinals of $\N$,
\item letting $\gamma=\sup_{i<\omega}\gg_i$, $\rho_\omega(\N)<\gg$, 
\item for some $k<\omega$, $\N$ is a sound $\Sigma$-mouse over $\N|\gg_k$, 
\item for every cutpoint cardinal $\eta<\gamma$, $\N|(\eta^+)^\N=\W(\N|\eta)$, and
\item letting $\Lambda$ be the $(\omega, \omega_1, \omega_1)$-strategy of $\N$, $Code(\Lambda) \not \in \Gamma$ and $\Lambda$ is $\Gamma$-fullness preserving (i.e., clause 4 holds for any $\Lambda$-iterate of $\N$). 
\end{enumerate}
\end{lemma}

We now begin the proof of \rlem{a mouse beyond n}. We continue with previous subsection's notation and start working in $\N^*_z$. Our aim is to use fully backgrounded constructions of $\N^*_z$ to produce a mouse $\N^*$ such that for some $l$, $\N=\mathcal{C}_l(\N^*)$ has the desired properties. Let $\la\eta_i : i<\omega\ra$ be the first $\omega$ cardinals below $\kappa$ such that for every $i<\omega$, $\W(\N^*_z|\eta_i)\models ``\eta_i$ is a Woodin cardinal" (it follows from \rlem{many woodin cardinal} that there are such cardinals). Let now $\la \N_i: i<\omega\ra$ be a sequence constructed according to the following rules:
\begin{enumerate}
\item $\N_0=(\mathcal{J}^{\vec{E}, \Sigma})^{\N^*_z|\eta_0}$,
\item $\N_{i+1}=(\mathcal{J}^{\vec{E}, \Sigma}[\N_i])^{\N^*_z|\eta_{i+1}}$.
\end{enumerate}
Let $\N_\omega=\cup_{i<\omega}\N_i$. \\

\textit{Claim 1.} For every $i<\omega$, $\N_\omega\models ``\eta_{i}$ is a Woodin cardinal" and $\N_{\omega}|(\eta_i^+)^{\N_\omega}=\W(\N_i)$.\\
\begin{proof}
It is enough to show that
\begin{enumerate}
\item $\N_{i+1}\models ``\eta_i$ is a Woodin cardinal",
\item no level of $\N_{i+1}$ projects across $\eta_i$, and
\item $\N_{i+1}|(\eta_i^+)^{\N_{i+1}}=\W(\N_i)$.
\end{enumerate}
To show 1-3, it is enough to show that if
$\Q\trianglelefteq \N_{i+1}$ is such that  $\rho_\omega(\Q)\leq\eta_i$ then the fragment of the iteration strategy of $\Q$ that acts on trees above $\eta_i$ is coded by a set in $\Gamma$ (this is simply because $\N_{i+1}$ is $\Gamma$-full). Fix then $i$ and let $\Q\trianglelefteq \N_{i+1}$ be such that  $\rho_\omega(\Q)\leq\eta_i$. Let $\xi$ be such that if $\S$ is the $\xi$th model of the fully backgrounded construction producing $\N_{i+1}$ then $\Q$ is the core of $\S$. Let $\pi: \Q\rightarrow \S$ be the uncollapse map. It is a fine structural map but that is irrelevant and we suppress this point. 

Let $\nu<\eta_{i+1}$ be a cardinal such that $\S$ is the $\xi$th model of the full background construction of $\N^*_z|\nu$. Let $\Psi$ be the fragment of $\Sigma_z$ that acts on non-dropping trees that are based on $\N^*_z|(\nu^+)^{\N^*_z}$ and are above $\eta_i$. We have that $\Psi$ induces an iteration strategy $\Psi^*$ for $\S$ and that $\pi$-pullback of $\Psi^*$ is an iteration strategy for $\Q$. It is then enough to show that $Code(\Psi)\in \Gamma$. 

Notice that whenever $\T$ is a tree on $\N^*_z|(\nu^+)^{\N^*_z}$ according to $\Psi$ and $b=\Psi(\T)$ then $\Q(b, \T)$ is defined. Also, notice that because of our choice of $\eta_{i+1}$, for any such $\T$ and $b$, $\Q(b, \T)\trianglelefteq \W(\M(\T))$. Because the function $a\rightarrow \W(a)$ is coded by a set in $\Gamma$, we have that $Code(\Psi)\in \Gamma$.  \\ 
\end{proof}

\textit{Claim 2.} There is $\Q\trianglelefteq (\mathcal{J}^{\vec{E}, \Sigma}(\N_\omega))^{\N^*_z}$ such that $\rho_\omega(\Q)<\eta_\omega$.\\
\begin{proof} To see this suppose not. Let $\R=(\mathcal{J}^{\vec{E}, \Sigma}(\N_\omega))^{\N^*_z}$. It follows from universality of $\R$ (with respect to $\Sigma$-mice that have iteration strategies in $\Gamma$), we have that 
\begin{center}
$Lp^{\Gamma_1, \Sigma}(\N_\omega)\trianglelefteq \R$.
\end{center}
It follows from our choice of $\Gamma_1$ and from our hypothesis that $Lp^{\Gamma_1, \Sigma}(\N_0)\trianglelefteq \N_\omega$. Notice that if $\M(a)$ is defined for some $a$ then because of our choice of $\Gamma_1$, $\M(a)\trianglelefteq Lp^{\Gamma_1, \Sigma}(a)$. 

We claim that $\M(\N_0)$ is defined. To see this, notice that $x$ is generic over $\mathcal{J}[\N_0]$ for the extender algebra at $\eta_0$. Hence, if $g\subseteq Coll(\omega, \eta_0)$ is $Lp^{\Gamma_1, \Sigma}(\N_0)$-generic such that $x\in Lp^{\Gamma_1, \Sigma}(\N_0)[g]$, then by the results of Section 2.11 of \cite{ATHM}, we have that 
\begin{center}
$Lp^{\Gamma_1, \Sigma}(\N_0)[g]=Lp^{\Gamma_1, \Sigma}(\N_0[g])$.
\end{center}
But now, because $x\in \mathcal{J}[\N_0][g]$, we have that $\M(\N_0[g])$ is defined, and by our choice of $\Gamma_1$, we have that $\M(\N_0[g])\trianglelefteq Lp^{\Gamma_1, \Sigma}(\N_0[g])$. Again using the results of Section 2.11 of \cite{ATHM}, we have that some initial segment of $Lp^{\Gamma_1, \Sigma}(\N_0)$ has an iteration strategy which is not coded by a set of reals in $\Gamma$. Hence, $\M(\N_0)$ is defined. 

Because $\M(\N_0)$ is defined, we have that $\M(\N_0)\trianglelefteq Lp^{\Gamma_1, \Sigma}(\N_0)$ and therefore, $\M(\N_0)\trianglelefteq \N_\omega$. However, it follows from the proof of Claim 1 that all initial segments of $\N_\omega$ projecting to $\eta_0$ have an iteration strategy coded by a set in $\Gamma$. This implies that $\M(\N_0)$ has an iteration strategy coded by a set in $\Gamma$, contradiction!
\end{proof}

Let now $\N^*\trianglelefteq Lp(\N_\omega)$ be least such that $\rho_\omega(\N^*)<\eta_\omega$. Let $l$ be least such that $\rho_l(\N^*)<\eta_\omega$ and let  $k$ be least such that $\rho_{l}(\N^*)<\eta_k$. In what follows, we will regard $\N^*$ as a $\Sigma$-mouse over $\N^*|\eta_k$. We let $\N=\mathbb{C}_{l}(\N^*)$. Thus, $\N$ is sound (as a $\Sigma$-mouse over $\N|\eta_k$). We let $\la \gamma_i : i<\omega\ra$ be the Woodin cardinals of $\N$ and $\gamma=\sup_{i<\omega}\gg_i$. Let $\Lambda$ be the $(\omega, \omega_1, \omega_1)$-strategy of $\N$ induced by $\Sigma_z$. Notice that $Code(\Lambda) \notin \Gamma$ because otherwise, since $\N_\omega$ is $\Gamma$-full, $\N \trianglelefteq \N_\omega \trianglelefteq \N^*$.\\

\textit{Claim 3.} $\Lambda$ is $\Gamma$-fullness preserving.\\
\begin{proof}
To see this fix $\N_1$ which is a $\Lambda$-iterate of $\N$ via $\VT$ such that the iteration embedding $i: \N\rightarrow \N_1$ exists. If $\N_1$ isn't $\Gamma$-full then there is a cutpoint $\nu$ of $\N_1$ and a sound $\Sigma$-mouse $\Q$ over $\N_1|\nu$ with $(\omega, \omega_1)$-iteration strategy $\Psi$ such that $Code(\Psi) \in \Gamma$, $\rho_\omega(\Q)=\nu$ and $\Q\ntrianglelefteq \N_1$.\\

\textit{Subclaim.} $\Psi$ can be extended to an $(\omega, \omega_1, \omega_1)$-iteration strategy.\\
\begin{proof}
We can find a good pointclass $\Gamma^*$ such that $Code(\Psi)\in \utilde{\Delta}_{\Gamma^*}$. Using Theorem 1.2.9 of \cite{ATHM}, we can find $(\N^*_y, \Sigma_y, \d_y)$ that Suslin captures $Code(\Psi)$. Notice that $\Sigma_y$ is an $(\omega, \omega_1, \omega_1)$-iteration strategy. It follows from universality that $\Q\trianglelefteq (\mathcal{J}^{\vec{E}, \Sigma}[\N_1|\nu])^{\N^*_y|\d_y}$. Hence, $\Q$ has an $(\omega, \omega_1, \omega_1)$-iteration strategy $\Psi^+$. Because $\Psi$ is the unique $(\omega, \omega_1)$-iteration strategy of $\Q$, we have that $\Psi^+$ extends $\Psi$. 
\end{proof}

We now compare $\Q$ with $\N_1$. Let $\S$ be the comparison tree on the $\Q$ side with last model $\Q^*$ and $\T$ be the comparison tree on the $\N_1$ side with last model $\N_1^*$. Because $\Q\ntrianglelefteq \N_1$, we must have that $\N_1^* \trianglelefteq \Q^*$ and $\pi^{\T}:\N_1\rightarrow \N_1^*$ exists. Because the $(\omega, \omega_1)$-fragment of $\Lambda$ is the unique $(\omega, \omega_1)$-iteration strategy of $\N$, we must have that it is the $\pi^\T\circ i$-pullback of $\Psi_{\N_1^*, \VT^\frown \S}$ (recall that this is the strategy of $\N_1^*$ induced by $\Psi$). This implies that $\Lambda \in \Gamma$,  contradiction.
\end{proof}

It is now clear that $(\N, \Lambda)$ is as desired. This completes the proof of \rlem{a mouse beyond n}.

\subsection{A Prikry generic}

In this subsection, while working in $N$, we define a Prikry forcing with the property that the generic object produces a sound countably iterable $\Sigma$-mouse $\R$ over $\bR$ such that $\R\in N$ and extends $(Lp^\Sigma(\bR))^N$. Clearly this is a contradiction.

We now start working in $N$. We now describe a function $f:\mathcal{D}^{<\omega}\rightarrow HC$ such that if $G\subseteq \mathbb{P}^f$ is $N$-generic then $f(G)$ is a $\Sigma$-premouse such that certain $\mathcal{J}^{\vec{E},\Sigma}$-construction of it is an initial segment of some $\Lambda$-iterate of $\N$. 

Following \cite{ATHM}, we say $\Q$ is $\Sigma$-\textit{suitable} (in $N$) if for some ordinal $\d$
\begin{enumerate}
\item $\d$ is the unique Woodin cardinal of $\Q$,
\item $o(\Q)=\sup_{n<\omega} (\d^{+n})^\Q$,
\item $\Q$ is full with respect to $\Sigma$-mice, i.e., for any cutpoint $\eta$, $Lp^\Sigma(\Q|\eta)\trianglelefteq \Q$.
\end{enumerate}
We let $\d^\Q$ be the Woodin cardinal of $\Q$. Similarly we can define the notion of a $\Sigma$-suitable $\Q$ over any set $a$. In particular, if $\Q$ is $\Sigma$-suitable and $\R$ is $\Sigma$-suitable over $\Q$ then $\R\models ``\d^\Q$ is a Woodin cardinal". Because we will only deal with $\Sigma$-suitable structures, we omit $\Sigma$ and just say suitable instead of $\Sigma$-suitable.  

A normal iteration tree $\U$ on a suitable $\P$ is \textit{short} if for all limit $\xi \leq lh(\U)$,
\begin{center}
 $Lp^\Sigma(\M(\U|\xi)) \vDash ``\delta(\U|\xi)$ is not Woodin".
\end{center}
 Otherwise, we say that $\U$ is \textit{maximal}. We say that a suitable $\P$ is \textit{short tree iterable} if for any short tree $\T$ on $\P$, there is a cofinal wellfounded branch $b$ such that $\Q(b,\T)$ exists and if $\pi^\T_b:\P\rightarrow \M^\T_b$ exists then $\M^\T_b$ is suitable.

Write $\P_y$ for the premouse coded by the real $y$. Let $a$ be countable transitive and $d\in \mathcal{D}$ be such that $a$ is coded by a real recursive in $d$. Put
\begin{center}
$\F^d_a$ = $\left\{ \P_z : \ z \leq_T d, \P_z \textrm{ is a short-tree iterable suitable premouse over } a\right\}$
\end{center}

\begin{lemma}\label{cone lemma} For any fixed $a$, there is a cone of $d$ such that  $\F^d_a \neq \emptyset$.
\end{lemma}
\begin{proof}
If not, the failure of the statement in the claim is a $\Sigma_1$ statement. Call this statement $\phi[a]$. Using $\utilde{\Sigma}^2_1(Code(\Sigma))$-reflection, we get a transitive model 
\begin{center}
$H \vDash ZF^- + \Theta = \Theta_{Code(\Sigma)} + \phi[a]$,
\end{center} 
 $\mathbb{R} \subseteq H$ and $\powerset(\bR)\cap H\subsetneq \utilde{\Delta}^2_1(Code(\Sigma))$. 

Let $\Gamma^*$ be a good pointclass beyond $H$. Such a $\Gamma^*$ exists by our assumption on $H$. We use Theorem 1.2.9 of \cite{ATHM} to get a triple $\langle N^*_w, \delta_w, \Sigma_w\rangle$ (for some real $w$) that Suslin captures the universal $\Gamma^*$ set. Using universality of fully backgrounded constructions and the proofs of the claims from the proof of \rlem{a mouse beyond n} (or the results of Section 3.2.2 of \cite{ATHM}), we conclude that the $(\mathcal{J}^{\vec{E}, \Sigma}[a])^{\N^*_w|\d_w}$ reaches a premouse $\Q_a$ such that in $H$, $\Q_a$ is short tree iterable and suitable (with respect to $H$). This contradicts our assumptions on $H$.
\end{proof}

For each $a$ and for each Turing degree $d$ from the cone of \rlem{cone lemma},  we can simultaneously compare all $\Q \in \F^d_a$ while doing the generic genericity iteration to make $d$ generic over the common part of the final model $\Q^{d,-}_a$. This process (hence $\Q^{d,-}_a$) depends only on $d$. Set
\begin{center}
$\Q^d_a = Lp_\omega^\Sigma(\Q^{d,-}_a)$ and $\delta^d_a = o(\Q^{d,-}_a)$.
\end{center}
Recall that we are working in $N$ (thus, we really have that $\Q^d_a = Lp_\omega^{\Gamma,\Sigma}(\Q^{d,-}_a)$).

\begin{lemma} The following statements are true (in $N$).
\end{lemma}
\begin{enumerate}
\item $\Q^d_a$ and $\delta^d_a$ depend only on $d$.
\item $\Q^{d,-}_a$ is $\Sigma$-full (no levels of $\Q^d_a$ project strictly below $\delta^d_a$).
\item $\Q^d_a \vDash \delta^d_a$ is Woodin.
\item $\powerset{(a)} \cap \Q^d_a = \powerset{(a)} \cap OD_\Sigma(a \cup \{a\})$ and $\powerset{(\delta^d_a)} \cap \Q^d_a = \powerset{(\delta^d_a)} \cap OD_\Sigma(Q^{d,-}_a \cup \{Q^{d,-}_a\})$.
\item $\delta^d_a = \omega_1^{L[S,d]}$.
\end{enumerate}
\begin{proof}
1-4 just follow from our definitions. We consider 5. Let $S$ be the tree of a $(\Sigma^2_1(Code(\Sigma)))^{N}$ scale on a universal $(\Sigma^2_1(Code(\Sigma)))^{N}$ set U. Suppose that in $L[S,d]$, the process producing $\Q^d_a$ stops at stage $\alpha < \omega_1^{L[S,d]}$. We then have that $\Q^d_a$ is countable in $L[S,d]$. The suitability of $\Q^d_a$ then implies that $\mathbb{R} \cap L[S,d] \subseteq \Q^d_a[d]$. It then follows that $\delta^d_a$, the Woodin of $\Q^d_a$, is countable in $L[S,d]$ while it is a cardinal in $\Q^d_a[d]$ (because the extender algebra of $\Q^d_a$ at $\d^d_a$ is $\delta^d_a$-cc). Hence, $\omega_1^{L[S, d]}$ is countable in $L[S, d]$, contradiction! 
\end{proof}

We now define $f:\mathcal{D}^{<\omega}\rightarrow HC$ by induction on $\mathcal{D}^n$. Fix $(\N, \Lambda)$ as in \rlem{a mouse beyond n} and let $k$ be as in clause 3. Below we use the notation of \rlem{a mouse beyond n}. We let $f(\emptyset)=\N|\gg_k$. Suppose we have defined $f\rest \mathcal{D}^{n+1}$. Given $p\in \mathcal{D}^{n+2}$, we let 
\begin{center}
$f(d)=\begin{cases}
\Q^d_{f(p\rest n+1)} &: f(p\rest n+1)\ \text{is countable in}\ L[d]\\
\emptyset &: \text{otherwise}
\end{cases}$
\end{center}

Suppose now that $G\subseteq \mathbb{P}^f$ is $N$-generic. Let $\Q_i=G^i$ and let $\Q_\omega=f(G)$. We let $\d_i$ be the largest Woodin cardinal of $\Q_i$. Without loss of generality, we assume that if $(\la d\ra, X)\in G$ then $\N$ is countable in $L[d]$.

Given an increasing function $h:\omega\rightarrow \omega$, we define $\la \Q^h_i, \Q^{h,*}_i : i<\omega\ra$ according to the following procedure:
\begin{enumerate}
\item $\Q^{h, *}_0$ is the output of $\mathcal{J}^{\vec{E}, \Sigma}[a]$ construction done in $\Q_{h(0)+1}$ using extenders with critical point $>\d_{h(0)}$.
\item $\Q^h_0=(Lp^\Sigma(\Q^{h, *}_0))^{\Q_{h(0)+2}}$.
\item $\Q^{h, *}_{ i+1}$ is the output of $\mathcal{J}^{\vec{E}, \Sigma}[\Q^{h}_i]$ construction done in $\Q_{h(i+1)+1}$ using extenders with critical point $>\d_{h(i+1)}$.
\item $\Q^{h}_{i+1}=(Lp^\Sigma(\Q^{h, *}_{ i+1 }))^{\Q_{h(i+1)+2}}$.
\end{enumerate}
We let $\Q^h_\omega=\cup_{i<\omega}\Q^h_i$.

\begin{lemma}\label{h lemma} For some increasing function $h:\omega\rightarrow\omega$ such that $h\in V$, $\Q^h_\omega$ is an initial segment of a $\Lambda$-iterate of $\N$.
\end{lemma}
\begin{proof} Let $\vec{d}=\la d_i: i<\omega\ra$ be the generic sequence of degrees given by $G$. We define $h$ recursively. It will have the property that $\Q^h_i$ is a $\Lambda$-iterate of $\N|(\gg_{k+i+1}^+)^\N$. While defining $h$, we also define a sequence $\vec{H}=\la \N_i, \U_i, b_i: i\in [-1, \omega)\ra$ such that
\begin{enumerate}
\item $\N_{-1}=\N$,
\item for each $i$, $\U_i$ is an iteration tree on $\N_i$ and $b_i=\Lambda(\oplus_{m<i+1} \U_m)$,
\item for each $i$, $\N_{i+1}=\M^{\U_i}_{b_i}$,
\item for each $i$, $\pi^{\U_i}_{b_i}$-exists and letting $\pi_{i, j}:\N_i\rightarrow \N_j$ be the composition of iteration embeddings, $\U_i$ is a tree based on $\N_i|[\pi_{-1, i}(\gg_{k+i}), \pi_{-1, i}(\gg_{k+i+1}))$,
\item $\Q^h_{i}=\N_i|(\pi_{-1, i}(\gg_{k+i+1}^+)^{\N})$
\item for each $i$, $h(i)=m+1$ where $m$ is the least integer such that $\vec{H}\rest i+1$ is countable in $L[d_m]$.
\end{enumerate} 
1-6 above tell us how to define the sequence. To see that we can always arrange 6, recall that $\vec{d}$ is cofinal in the set of degrees. To see that $h \in V$, recall that Prikry property implies that $\mathbb{P}^f$ doesn't add new reals. To see 5, notice that by our construction, $\vec{H}\rest i$ is generic over $\Q_i$ for the extender algebra at $\d_i$. 
\end{proof}

We let $h$ be as in \rlem{h lemma}. We let $\S_i=\Q_i^h$ and $\S_\omega=\Q_\omega^h$. Also, let $\S$ be the $\Lambda$-iterate of $\N$ such that $\S_\omega\trianglelefteq \S$. Because $\rho(\N)\leq \gg_k$, we have that $\rho(\S)\leq \gg_k$. Let $\la \eta_n: n<\omega\ra$ be the Woodin cardinals of $\S$. Let $\eta_\omega=\sup_{n<\omega}\eta_n$. Notice that in $V[G]$, $\S_\omega$ is $(\omega, \omega_1)$-iterable for short trees.

\indent We now have that there is $g\subseteq Col(\omega, <\eta_\omega)$-generic over $\S$ such that 
 \begin{center}
 $\cup_{n<\omega} \mathbb{R}^{\S[g\cap \eta_n]}=\mathbb{R}$. 
 \end{center} 
Next we perform an $S$-construction (see Section 2.11 of \cite{ATHM}, \cite{Selfiterability}, \cite{steel03} or \cite{steel08}) to translate $\mathcal{S}$ to a $\Sigma$-mouse over  $\mathbb{R}$. To see that the translation procedure works, let $\l=\Theta^{\mathcal{J}_\a(U, \bR)}$.  Notice that $\mathbb{P}^f\in \mathcal{J}_\l(\mathbb{R})$ and that all extenders of $\S$ above $\eta_\omega$ have critical point $>\l$. Thus, we can translate $\Sigma$-premice over $\mathcal{J}_\l[\S_\omega]$ to $\Sigma$-premice over $\mathcal{J}_\l(\mathbb{R})$. Let then $\W$ be the $\Sigma$-premouse over $\bR$ that is the result of translating $\S$ into a $\Sigma$-premouse over $\R$. 

\begin{lemma} $(Lp^\Sigma(\mathbb{R}))^N\inseg \W$.
\end{lemma}
\begin{proof}
Suppose $\W\trianglelefteq (Lp^\Sigma(\mathbb{R}))^N$. Notice that $\W$ is $OD^N_\Sigma$. Notice that $\N$ is the $\d_k$-core of $\S$. Let then $\tau$ be a name for a sound $\Sigma$-premouse over $\N|\d_k$ projecting to $\d_k$ such that it is always realized as the $\d_k$-core of the translation of $\W$ into an extension of $\Q_\omega^h$. Then $\tau$ is $OD_\Sigma^N$ and hence, there is $OD_\Sigma^N$ condition $(\emptyset, X)$ that decides $\tau$. It then follows that if $\N^*$ is the premouse given by $\tau$ and $(\emptyset, X)$ then $\N^*=\N$. But this implies that $\N\in OD_{\N|\d_k, \Sigma}^N$ and hence, by $N$-fullnes of $\N$, $\N\in \N$, contradiction.
\end{proof}

Let then $\mathcal{R}\insegeq \W$ be the first level of $\W$ such that 
\begin{center}
$(Lp^\Sigma(\mathbb{R}))^N\trianglelefteq \mathcal{R}$ and $\rho_\omega(\R)=\mathbb{R}$.
\end{center} 

The next two lemmas finishes the prove of \rthm{main theorem}.

\begin{lemma} $\mathcal{R}\in V$.
\end{lemma}
\begin{proof}
Suppose not. Using $DC$ we can find $\pi: H\rightarrow \mathcal{J}_{\mu}(\powerset(\mathbb{R}))$ such that $\mu>\Theta$ and $H$ is countable. We can further assume that $\Sigma, \Lambda, N\in rng(\pi)$. Let then $\bar{N}=\pi^{-1}(N)$. Let $g\subseteq Coll(\omega, \bar{N})$ be $H$-generic and let $g_1, g_2\subseteq \pi^{-1}(\mathbb{P}^f)\in H[g]$ be two different $\bar{N}$-generics. Let $\R_1$ and $\R_2$ be the versions of $\R$ defined for $\bar{N}$ using $g_1$ and $g_2$ respectively. Because both are $(\omega, \omega_1)$-iterable, we have that $\R_1=\R_2$. Hence, the version of $\R$ for $\bar{N}$ is $OD_{\Sigma^H}^{H[g]}$ and hence, it is in $H$.  
\end{proof}
It remains to show that $\mathcal{R}$ is in $N$ and countably iterable in $N$. Granted this, we obtain the desired contradiction, hence complete the proof of \rthm{main theorem}.

\begin{lemma}
$\mathcal{R} \in N$ and $\mathcal{R}$ is countably iterable in N.
\end{lemma}
\begin{proof}
First, we show $\mathcal{R} \in N$. We can assume $o(\mathcal{R})$ is limit and $\rho_1(\mathcal{R}) = \mathbb{R}$ (if not, look at the mastercode structure of $\mathcal{R}$). In V, we can write $\mathcal{R}   =  \cup_{\xi<o(\mathcal{R})} Th^{\mathcal{R}|\xi}(\mathbb{R})$. Notice that for all $\xi < o(\mathcal{R})$, $|Th^{\mathcal{R}|\xi}(\mathbb{R})|_w < |A|_w$ (where $A$ is the least $OD_{\Sigma}^N$ set of reals such that $A\not \in Lp^{\Sigma}(\bR))^N$). Since $\mathcal{R}$ is a well-ordered union of sets Wadge reducible to A, it follows from a theorem of Kechris that $\mathcal{R}$ is projective in A. This implies that $\mathcal{R} \in N$.
\\
\indent It remains to show $\mathcal{R}$ is countably iterable in $N$. Working in $N$, given $\sigma\in \powerset_{\omega_1}(\mathbb{R})$ we say $\sigma$ is bad if there is a non-iterable sound $\Sigma$-premouse $\W$ over $\sigma$ projecting to $\sigma$ and an embedding $\pi:\W\rightarrow \R$. Notice that (in V) $\R$ is countably $(\omega, \omega_1)$-iterable and hence, for each $\sigma\in \powerset_{\omega_1}(\mathbb{R})$ there is at most one such $\W$. We denote it by $\W(\sigma)$. 

To show that $\R$ is countably iterable in $N$, it is enough to show that for stationary many $\sigma$, $\W(\sigma)$ is undefined. Towards a contradiction assume that for a club $C$ of $\sigma$, $\W(\sigma)$ is defined. Then the set
\begin{center}
 $B=\{(\sigma, \W(\sigma)): \sigma\in C\}$.
 \end{center}
  is $OD_{\Sigma, u}^N$ for some real $u$. It follows that for every $\sigma\in C$ such that $u\in \sigma$, $\W(\sigma)\in (Lp^\Sigma(\sigma))^N$. Because for every $\sigma$, $\W(\sigma)$ has an $(\omega, \omega_1)$-iteration strategy in $V$, we get that $\W(\sigma)\trianglelefteq (Lp^{\Sigma}(\sigma))^{N}$, which is a contradiction.
\end{proof}

\bibliographystyle{plain}
\bibliography{Rmice}
\end{document}